# SEQUENTIAL CHANGE DETECTION REVISITED

By George V. Moustakides

*University of Patras*

In sequential change detection, existing performance measures differ significantly in the way they treat the time of change. By modeling this quantity as a random time, we introduce a general framework capable of capturing and better understanding most well-known criteria and also propose new ones. For a specific new criterion that constitutes an extension to Lorden's performance measure, we offer the optimum structure for detecting a change in the constant drift of a Brownian motion and a formula for the corresponding optimum performance.

**1. Introduction.** Suppose we are observing sequentially a process $\{\xi_t\}_{t>0}$, which up to and including time $\tau \geq 0$ follows the probability measure $\mathbb{P}_\infty$ and after $\tau$ it switches to an alternative regime $\mathbb{P}_0$. Parameter $\tau$ is the change-time and denotes the *last* time instant the process is under the nominal regime $\mathbb{P}_\infty$. The goal is to detect the change of measures as soon as possible, using a sequential scheme.

Any sequential test can be modeled as a stopping time (s.t.) $T$ adapted to the filtration $\{\mathcal{F}_t\}_{t\geq 0}$, where $\mathcal{F}_t = \sigma\{\xi_s, 0 < s \leq t\}$ for $t > 0$; and $\mathcal{F}_0$ is the trivial $\sigma$-algebra. We note that the process $\{\xi_t\}$ becomes available for $t > 0$ while the change-time $\tau$ can take upon the value 0 as well. This is because with $\tau = 0$ we would like to capture the case where all observations are under the alternative regime, whereas $\tau = \infty$ refers to the case where all observations are under the nominal regime. More generally, $\mathbb{P}_\tau$ denotes the probability measure induced by the change occurring at $\tau$ and $\mathbb{E}_\tau[\cdot]$ the corresponding expectation. In particular, if $X$ is an $\mathcal{F}_\infty$-measurable random variable and $\tau = t$ a *deterministic* time of change, then we can write

$$(1.1) \qquad \mathbb{E}_t[X] = \mathbb{E}_\infty[\mathbb{E}_0[X|\mathcal{F}_t]].$$

In developing optimum change detection algorithms, the first step consists in defining a suitable performance measure. Existing criteria basically quantify the detection delay $(T - \tau)^+$, where $x^+ = \max\{x, 0\}$, by considering alternative versions of its average. These definitions play an important role in the underlying mathematical model for the change-time $\tau$.

Currently we distinguish two major models for change-time. The first, introduced by Shiryayev (1978), assumes that $\tau$ is random with a known (exponential) prior. We can accompany this change-time model with the following performance measure $\mathcal{J}_S(T) = \mathbb{E}_\tau[T - \tau | T > \tau]$, that is, the average detection delay conditioned on the event that we stop after the change. Alternatively, we can consider $\tau$ to be deterministic and unknown and follow a worst-case scenario. There exist two possibilities. The first, proposed by Lorden (1971), considers the worst average delay $\mathcal{J}_L(T) = \sup_{0 \leq \tau} \operatorname{esssup} \mathbb{E}_\tau[(T - \tau)^+ | \mathcal{F}_\tau]$ conditioned on the least-favorable observations before the change. The second, due to Pollak (1985), uses the worst average delay $\mathcal{J}_P(T) = \sup_{0 \leq \tau} \mathbb{E}_\tau[T - \tau | T > \tau]$ conditioned on the event that we stop after $\tau$.

Shiryayev's Bayesian approach presents definite analytical advantages and has been the favorite underlying model in several existing optimality results as Poor (1998), Beibel (1996), Peskir and Shiryayev (2002), Karatzas (2003), Bayraktar and Dayanik (2006) and Bayraktar, Dayanik and Karatzas (2006). The two deterministic approaches on the other hand, although more analytically involved, are clearly more tractable from a practical point of view since they do not make any limiting assumptions.

As it will become evident in Section 3, the three performance criteria can be ordered as follows: $\mathcal{J}_S(T) \leq \mathcal{J}_P(T) \leq \mathcal{J}_L(T)$. Because of this property, there exist strong arguments against Lorden's measure as being overly pessimistic. Such claims, however, tend to be inconsistent with the fact that $\mathcal{J}_L(T)$, whenever it can be optimized, it gives rise to the CUSUM s.t., one of the most widely used change detection schemes in practice. Despite their similarity, Pollak's $\mathcal{J}_P(T)$ and Lorden's $\mathcal{J}_L(T)$ measure, as we are going to see, differ in a very essential way. In fact $\mathcal{J}_P(T)$, although not obvious at this point, will be shown to be closer to Shiryayev's $\mathcal{J}_S(T)$ measure than to Lorden's $\mathcal{J}_L(T)$.

In the next section we present a general approach for modeling the change-time $\tau$. The three measures presented previously will turn out to be special cases of our general setting corresponding to different levels and forms of prior knowledge. The understanding of their differences will give rise to a discussion concerning the suitability of each measure for the problem of interest and will explain, we believe in a convincing way, why Lorden's criterion, although seemingly more pessimistic than the other two, is more appropriate for the majority of change detection problems. Finally, we are going to introduce an additional criterion that constitutes an extension to Lorden's $\mathcal{J}_L(T)$ performance measure. For this case, we will also provide the optimum test for detecting a change in the constant drift of a Brownian motion and a formula for the corresponding optimum performance.

**2. A randomized change-time.** Suppose that nature, at every time instant $t$, consults the available information $\mathcal{F}_t$ and with some *randomization* probability decides whether it should continue using the nominal probability measure or switch to the alternative one. Consequently, let $\pi_t$ denote the randomization probability that there is a change at $t$ conditioned on the available information up to time $t$, that is $\pi_t = \mathbb{P}[\tau = t | \mathcal{F}_t]$. Clearly, $\pi_t$ is nonnegative and the process $\{\pi_t\}$ is $\{\mathcal{F}_t\}$-adapted.

We recall that time $\tau$ is usually considered in the literature as the *first* time instant under the alternative regime. With the current setting this is no longer possible. Indeed, since there is a decision involved whether to change the statistics or not, this decision must be made *before* any data under the alternative regime are produced. Therefore, $\tau$ denotes the time we *stop* using the nominal regime.

Consider now a process $\{\mathcal{X}_t\}_{t \geq 0}$, where $\mathcal{X}_t$ is nonnegative and $\mathcal{F}_\infty$-measurable (the process in not necessarily $\{\mathcal{F}_t\}$-adapted). We would like to compute the expectation of the random variable $\mathcal{X}_\tau$ which is the $\tau$-randomly-stopped version of $\{\mathcal{X}_t\}$, but we are interested only in *finite* values of $\tau$. In other words we would like to find $\mathbb{E}_\tau[\mathcal{X}_\tau | \tau < \infty]$. Using (1.1) and that $\pi_t$ is $\mathcal{F}_t$-measurable, we can write

$$\mathbb{E}_\tau[\mathcal{X}_\tau \mathbb{1}_{\{\tau < \infty\}}] = \sum_{t=0}^{\infty} \mathbb{E}_t[\mathcal{X}_t \pi_t] = \sum_{t=0}^{\infty} \mathbb{E}_\infty[\mathbb{E}_0[\mathcal{X}_t | \mathcal{F}_t] \pi_t].$$

Substituting $\mathcal{X}_t = 1$ in the previous relation, we obtain $\mathbb{P}_\tau[\tau < \infty] = \sum_{t=0}^{\infty} \mathbb{E}_\infty[\pi_t]$, which is an expression for the probability of stopping at finite time. Combining the two outcomes leads to

$$\mathbb{E}_\tau[\mathcal{X}_\tau | \tau < \infty] = \frac{\sum_{t=0}^{\infty} \mathbb{E}_\infty[\mathbb{E}_0[\mathcal{X}_t | \mathcal{F}_t] \pi_t]}{\sum_{t=0}^{\infty} \mathbb{E}_\infty[\pi_t]}.$$

From now on, and without loss of generality, we make the simplifying assumption that $\mathbb{P}_\tau[\tau < \infty] = 1$ (otherwise divide each $\pi_t$ with $\mathbb{P}_\tau[\tau < \infty]$). Under this assumption, we have

$$(2.1) \qquad \mathbb{E}_\tau[\mathcal{X}_\tau] = \sum_{t=0}^{\infty} \mathbb{E}_\infty[\mathbb{E}_0[\mathcal{X}_t | \mathcal{F}_t] \pi_t].$$

Let us summarize our change-time model. We are given a time increasing information (filtration) $\{\mathcal{F}_t\}_{t \geq 0}$ with $\mathcal{F}_0$ being the trivial $\sigma$-algebra, and a sequence of $\{\mathcal{F}_t\}$-adapted probabilities $\{\pi_t\}$. Quantity $\pi_t$ denotes the history dependent randomization probability that $t$ is the last time instant we obtain information under the nominal probability $\mathbb{P}_\infty$ and at the next time instant the new information will follow the alternative measure $\mathbb{P}_0$. For a process $\{\mathcal{X}_t\}$ with $\mathcal{X}_t$ being nonnegative and $\mathcal{F}_\infty$-measurable, we define the expectation of the $\tau$-randomly-stopped process $\mathcal{X}_\tau$ with respect to the measure induced by the change, with the help of (2.1).

2.1. *Decomposition of the change-time statistics.* The process $\{\pi_t\}$ can be decomposed as $\pi_t = \varpi_t p_t$ where $\{\varpi_t\}$ is a deterministic sequence of probabilities defined as

$$(2.2) \qquad \varpi_t = \mathbb{E}_\infty[\pi_t] \qquad \text{therefore} \sum_{t=0}^\infty \varpi_t = \mathbb{P}_\tau[\tau < \infty] = 1;$$

and $\{p_t\}$ a nonnegative $\{\mathcal{F}_t\}$-adapted process defined for $\varpi_t > 0$ as

$$(2.3) \qquad p_t = \frac{\pi_t}{\varpi_t} \qquad \text{therefore} \ \mathbb{E}_\infty[p_t] = 1,$$

while for $\varpi_t = 0$, we can arbitrarily set $p_t = 1$. Quantity $\varpi_t$ expresses the *aggregate* probability that $\tau$ will stop at $t$, whereas $p_t$ describes how this probability is *distributed* among the possible events that can occur up to time $t$. Since $\mathcal{F}_0$ is the trivial $\sigma$-algebra, $\pi_0$ is deterministic, therefore $\varpi_0 = \pi_0$ and $p_0 = 1$. Clearly $\varpi_0$ expresses the probability that the change takes place before the statistician obtains any information.

**3. Performance measure and optimization criterion.** If $T$ is an $\{\mathcal{F}_t\}$-adapted s.t. used by the statistician to detect the change, then we are interested in defining a measure that quantifies its performance. Following the idea of Lorden (1971) and Pollak (1985), we propose the use of

$$\mathcal{J}(T) = \mathbb{E}_\tau[T - \tau | T > \tau],$$

namely, the average detection delay conditioned on the event that we stop after $\tau$. Of course this measure makes sense for *finite* values of $\tau$ because a change at infinity is regarded as "no change." Since $(T - t)^+$ and $\mathbb{1}_{\{T>t\}}$ are nonnegative and $\mathcal{F}_\infty$-measurable, by using (2.1) our measure can be written as

$$(3.1) \qquad \begin{aligned} \mathcal{J}(T) &= \frac{\sum_{t=0}^\infty \mathbb{E}_\infty[\pi_t \mathbb{E}_0[(T-t)^+ | \mathcal{F}_t]]}{\sum_{t=0}^\infty \mathbb{E}_\infty[\pi_t \mathbb{1}_{\{T>t\}}]} \\ &= \frac{\sum_{t=0}^\infty \varpi_t \mathbb{E}_\infty[p_t \mathbb{E}_0[(T-t)^+ | \mathcal{F}_t]]}{\sum_{t=0}^\infty \varpi_t \mathbb{E}_\infty[p_t \mathbb{1}_{\{T>t\}}]}. \end{aligned}$$

If we are interested in finding an optimum $T$, then we must minimize $\mathcal{J}(T)$ with respect to $T$, controlling at the same time the rate of false alarms. Similarly to Lorden (1971) and Pollak (1985), we propose the following constrained optimization with respect to $T$:

$$(3.2) \qquad \inf_T \mathcal{J}(T) \quad \text{subject to} \quad \mathbb{E}_\infty[T] \geq \gamma.$$

In other words, we minimize the conditional average detection delay, subject to the constraint that the average period between false alarms is no less than a given value $\gamma \geq 0$. The performance measure, as we can see from (3.1), requires complete knowledge of the two processes $\{\varpi_t\}$ and $\{p_t\}$. In the next subsection we extend our definition to include cases where the statistics of $\tau$ are not exactly known or they are limited to special cases.

3.1. *Special cases and uncertainty classes.* If $\{\varpi_t\}, \{p_t\}$ are not known exactly and instead we have available an uncertainty class $\mathcal{T}$ for $\tau$, then we can extend the definition of our performance measure by adopting a worst-case approach of the form $\sup_{\tau \in \mathcal{T}} \mathcal{J}(T)$, while (3.2) can be replaced by the following min-max constrained optimization problem:

$$(3.3) \qquad \inf_T \sup_{\tau \in \mathcal{T}} \mathcal{J}(T) \quad \text{subject to} \quad \mathbb{E}_\infty[T] \geq \gamma.$$

Next, we are going to identify the particular form of our criterion for specific change-time classes. In order to facilitate our presentation, we first introduce a technical lemma.

LEMMA 1. *Let $\{\varpi_t\}$ and $\{p_t\}$ be the processes defined in Section 2.1 satisfying (2.2) and (2.3), respectively. If $\{a_t\}$, $\{b_t\}$ are two nonnegative deterministic sequences then*

$$(3.4) \qquad \sup_{\{\varpi_t\}} \frac{\sum_{t=0}^\infty \varpi_t a_t}{\sum_{t=0}^\infty \varpi_t b_t} = \sup_{0 \leq t} \frac{a_t}{b_t},$$

*where, for $a_t = b_t = 0$ we define the ratio $a_t/b_t = 0$. Furthermore, if $x_t, y_t$ are two nonnegative and $\mathcal{F}_t$-measurable random variables then*

$$(3.5) \qquad \sup_{p_t} \frac{\mathbb{E}_\infty[p_t x_t]}{\mathbb{E}_\infty[p_t y_t]} = \operatorname{esssup} \frac{x_t}{y_t},$$

*where, as before, when $x_t = y_t = 0$ we define the ratio $x_t/y_t = 0$.*

PROOF. To prove (3.4) notice that since $a_t \leq \{\sup_{0 \leq t}(a_t/b_t)\}b_t$ we conclude that for any sequence $\{\varpi_t\}$ we have

$$(3.6) \qquad \frac{\sum_{t=0}^\infty \varpi_t a_t}{\sum_{t=0}^\infty \varpi_t b_t} \leq \sup_{0 \leq t} \frac{a_t}{b_t}.$$

The upper bound in (3.6) is attainable by a sequence $\{\varpi_t\}$ that places all its probability mass on the time instant(s) that attain the supremum. If the supremum is attained in the limit, then for every $\epsilon > 0$ we can find a sequence $\{\varpi_t\}$ that depends on $\epsilon$, such that the left-hand side in (3.6) is $\epsilon$ close to the right-hand side.

Similar arguments apply for the proof of (3.5). Notice that, for every $p_t \geq 0$ satisfying $\mathbb{E}_\infty[p_t] = 1$ the combination $\mathbb{E}_\infty[p_t \cdot]$ defines a probability measure on $\mathcal{F}_t$ which is absolutely continuous with respect to $\mathbb{P}_\infty$. Since $x_t \leq \{\text{essup}(x_t/y_t)\}y_t$, $\mathbb{P}_\infty$-a.s., this leads to

$$\frac{\mathbb{E}_\infty[p_t x_t]}{\mathbb{E}_\infty[p_t y_t]} \leq \text{essup}\,\frac{x_t}{y_t}.$$

The upper bound is attainable by a probability measure $\mathbb{E}_\infty[p_t \cdot]$ that places all its mass on the event(s) that attain the essup, or we use limiting arguments if the essup is attained as a limit. $\square$

Let us now proceed with the presentation of specific special cases and uncertainty classes regarding the two processes $\{\varpi_t\}, \{p_t\}$.

*Case of known $\varpi_t$ and $p_t = 1$.* Here, by selecting $p_t = 1$, we limit our general change-time model to the case where the probability that the change will occur at $t$ is *independent from the observed history* $\mathcal{F}_t$. The corresponding performance measure simplifies to the following expression

$$(3.7) \qquad \mathcal{J}_\mathrm{S}(T) = \mathcal{J}(T)|_{p_t=1} = \frac{\sum_{t=0}^\infty \varpi_t \mathbb{E}_t[(T-t)^+]}{\sum_{t=0}^\infty \varpi_t \mathbb{P}_\infty[T > t]},$$

where we used (1.1) to replace $\mathbb{E}_\infty[\mathbb{E}_0[(T-t)^+|\mathcal{F}_t]]$ with $\mathbb{E}_t[(T-t)^+]$. There is no uncertainty class involved, we have simply limited the change-time $\tau$ to this special case. We recall that Shiryaev (1978) first introduced this model for the particular selection $\varpi_t = (1-\delta)\delta^t$.

*Case of arbitrary $\varpi_t$ and $p_t = 1$.* We continue using the same model of the previous case, but now we let $\{\varpi_t\}$ be an arbitrary sequence of probabilities satisfying, according to (2.2), $\sum_{t=0}^\infty \varpi_t = 1$. Using (3.4) from Lemma 1 and (1.1), it is straightforward to prove that

$$(3.8) \qquad \mathcal{J}_\mathrm{P}(T) = \sup_{\{\varpi_t\}} \mathcal{J}(T)|_{p_t=1} = \sup_{0 \leq t} \mathbb{E}_t[T - t|T > t].$$

By considering arbitrary $\{\varpi_t\}$, we recover Pollak's performance measure. From the way $\mathcal{J}_\mathrm{P}(T)$ is defined, it is evident that $\mathcal{J}_\mathrm{S}(T) \leq \mathcal{J}_\mathrm{P}(T)$.

Regarding the minimization of $\mathcal{J}_\mathrm{P}(T)$ with respect to the s.t. $T$, Pollak (1985) proposed the solution of the constrained optimization problem in (3.3). As candidate optimum s.t. for i.i.d. observations he suggested the Shiryaev–Roberts stopping rule. Pollak was able to demonstrate asymptotic optimality (as $\gamma \to \infty$) for this test. Regarding nonasymptotic optimality of the Shiryayev–Robert s.t. with respect to this criterion, see Mei (2006).

*Case of arbitrary $\varpi_t$ and arbitrary $p_t$.* Here the probability to stop at time $t$ depends on the observed history $\mathcal{F}_t$, we thus return to our general change-time model, but we assume complete lack of knowledge for the change time probabilities. In order to find the worst-case performance, we need to maximize $\mathcal{J}(T)$ with respect to both processes $\{\varpi_t\}$ and $\{p_t\}$. We have the following lemma that treats this problem.

LEMMA 2. *Let $\{\varpi_t\}$ and $\{p_t\}$ be defined as in Section* 2.1 *satisfying* (2.2) *and* (2.3) *respectively, then*

$$(3.9) \qquad \mathcal{J}_{\mathrm{L}}(T) = \sup_{\{\varpi_t\},\{p_t\}} \mathcal{J}(T) = \sup_{0 \leq t} \operatorname{ess\,sup} \mathbb{E}_t[(T-t)^+|\mathcal{F}_t].$$

PROOF. Using (3.4) from Lemma 1, for any given sequence $\{p_t\}$ we have

$$\sup_{\{\varpi_t\}} \mathcal{J}(T) = \sup_{0 \leq t} \frac{\mathbb{E}_\infty[p_t \mathbb{E}_0[(T-t)^+|\mathcal{F}_t]]}{\mathbb{E}_\infty[p_t \mathbb{1}_{\{T>t\}}]}.$$

Using the fact that we can change the order of two consecutive maximizations, we have

$$\begin{aligned}
\sup_{\{p_t\},\{\varpi_t\}} \mathcal{J}(T) &= \sup_{\{p_t\}} \sup_{0 \leq t} \frac{\mathbb{E}_\infty[p_t \mathbb{E}_0[(T-t)^+|\mathcal{F}_t]]}{\mathbb{E}_\infty[p_t \mathbb{1}_{\{T>t\}}]} \\
&= \sup_{0 \leq t} \sup_{p_t} \frac{\mathbb{E}_\infty[p_t \mathbb{E}_0[(T-t)^+|\mathcal{F}_t]]}{\mathbb{E}_\infty[p_t \mathbb{1}_{\{T>t\}}]} \\
&= \sup_{0 \leq t} \operatorname{ess\,sup} \mathbb{E}_0[(T-t)^+|\mathcal{F}_t] \\
&= \sup_{0 \leq t} \operatorname{ess\,sup} \mathbb{E}_t[(T-t)^+|\mathcal{F}_t],
\end{aligned}$$

where for the third equality we used (3.5) from Lemma 1 and for the last equality the fact that $\mathbb{E}_0[\cdot|\mathcal{F}_t] = \mathbb{E}_t[\cdot|\mathcal{F}_t]$. This concludes the proof. □

Here we recover Lorden's performance measure. It is clear that $\mathcal{J}_{\mathrm{P}}(T) \leq \mathcal{J}_{\mathrm{L}}(T)$, since for Lorden's measure we maximize over $\{p_t\}$ while in $\mathcal{J}_{\mathrm{P}}(T)$ we consider $p_t = 1$. As it was demonstrated in Moustakides (1986) and Ritov (1990), solving the optimization problem in (3.3) for Lorden's criterion and for i.i.d. observations, gives rise to the CUSUM test proposed by Page (1954). It is interesting to mention that Ritov (1990) based his proof of optimality on a change-time formulation, similar to the one proposed here.

A slight variation of the previous uncertainty class consists in assuming that the change cannot occur outside a sequence $\{t_n\}_{n \geq 0}$ of known time instants. In other words, we have $\varpi_t = 0$ if $t \notin \{t_0, t_1, \dots\}$. This modifies the previous criterion in the following way

$$(3.10) \qquad \mathcal{J}_{\mathrm{EL}}(T) = \sup_{\{\varpi_{t_n}\}, \{p_{t_n}\}} \mathcal{J}(T) = \sup_{0 \le n} \mathrm{essup}\, \mathbb{E}_{t_n}[(T - t_n)^+ | \mathcal{F}_{t_n}].$$

With a more accurate description of the time instants where the change can occur, one might expect to improve detection as compared to the CUSUM test. This measure is presented for the first time and will be treated in detail and under a more interesting frame in Section 4.

It is also possible to examine, under the general model, the case where $\{\varpi_t\}$ is known and $\{p_t\}$ unknown or, alternatively, $\{\varpi_t\}$ unknown and $\{p_t\}$ known. Clearly, the first case could be regarded as an extension of Shiryaev's approach to the general change-time model proposed here. Unfortunately both cases lead to rather complicated performance criteria, we, therefore, omit the corresponding analysis.

*Discussion.* From the preceding presentation it is evident that the three performance measures are ordered in the following way:

$$\mathcal{J}_{\mathrm{S}}(T) \le \mathcal{J}_{\mathrm{P}}(T) \le \mathcal{J}_{\mathrm{L}}(T),$$

giving the impression that Lorden's criterion is more pessimistic than Shiryaev's and Pollak's. This conclusion, however, is misleading since the underlying change-time model for the Shiryaev and Pollak criterion is completely different and significantly more limited than Lorden's. We recall that $\mathcal{J}_{\mathrm{S}}(T)$ and $\mathcal{J}_{\mathrm{P}}(T)$ rely on the assumption that the change at time $t$ is triggered with a probability that *does not depend* on the observed history $\mathcal{F}_t$. In practice there are clearly applications where this assumption is false and where it is more realistic to assume that the observations supply at least some partial information about the events that can trigger the change. Therefore, whenever we adopt this logic, Lorden's performance measure becomes more suitable than Shiryaev's and Pollak's. The same way $\mathcal{J}_{\mathrm{P}}(T)$ is preferable to $\mathcal{J}_{\mathrm{S}}(T)$ when there is no prior knowledge of $\{\varpi_t\}$ [despite the fact that $\mathcal{J}_{\mathrm{P}}(T)$ is more "pessimistic" than $\mathcal{J}_{\mathrm{S}}(T)$], we can also argue that $\mathcal{J}_{\mathrm{L}}(T)$ is preferable to $\mathcal{J}_{\mathrm{S}}(T)$ and $\mathcal{J}_{\mathrm{P}}(T)$ for problems where we need to follow the general change-time model and there is no prior knowledge regarding the change-point mechanism. Even if we still insist that $\mathcal{J}_{\mathrm{L}}(T)$ is overly pessimistic, it has now become clear that $\mathcal{J}_{\mathrm{S}}(T)$ and $\mathcal{J}_{\mathrm{P}}(T)$ are not the right alternatives, since they correspond to a drastically different change-time model.

Our previous arguments also suggest a word of caution when evaluating or comparing performances through Monte Carlo simulations. Selecting the time of change in an arbitrary way that has no relation with the observation sequence, is equivalent to adopting the restrictive change-time model with $p_t = 1$. This in turn is expected to favor tests that rely on this specific selection.

**4. Change at observable random times.** Let us now attempt a different parametrization of the change-time $\tau$. Suppose that in addition to the process $\{\xi_t\}_{t>0}$ we also observe a strictly increasing sequence of random times $\{\tau_n\}, n = 1, 2, \ldots$. These times correspond to occurrences of random events that *can trigger the change* of measures. In other words, we make the assumption that the change can occur *only* at the observable time instants $\{\tau_n\}$. We would like to emphasize that we consider the flow of the observation sequence $\{\xi_t\}$ to be continuous and not synchronized in any sense with the random times $\{\tau_n\}$. It is, therefore, clear that detection can be performed at *any* time instant, that is, even between occurrences.

There are interesting applications that can be modeled with this setup. For example, earthquake damage detection in structures, where earthquakes occurring at (observable) random times can trigger a change (damage), while detection is performed by continuously acquiring vibration measurements from the structure. Similar application is the detection of a change in financial data after "major importance events" or, as reported by Rodionov and Overland (2005), detection of regime shifts in sea ecosystems due to (observable) changes in the climate system.

Let us now relate our problem to the change-time model introduced in the previous section. Consider the strictly increasing sequence of occurrence times $\{\tau_n\}$, $n = 1, 2, \ldots$. Since we assume that observations are available *after* time 0, it is clear that $\tau_1 > 0$, therefore, we arbitrarily include $\tau_0 = 0$ into our sequence. Notice that $\tau_0$ does not necessarily correspond to a real occurrence. This term is needed to account for the case where the change took place before any observation was taken. If $\mathcal{N}_t$ denotes the number of observed occurrences up to (and including) time $t$, that is,

$$\mathcal{N}_t = \sup_{0 \le n}\{n : \tau_n \le t\},$$

then we can define our filtration $\{\mathcal{F}_t\}$ as $\mathcal{F}_t = \sigma\{\xi_s, \mathcal{N}_s, \ 0 < s \le t\}$ and $\mathcal{F}_0$ to be the trivial $\sigma$-algebra. With this filtration the random times $\tau_n$ are transformed immediately into s.t. adapted to $\{\mathcal{F}_t\}$ (since by consulting the history $\mathcal{F}_t$ we can directly deduce whether $\tau_n \le t$ is true).

The probability $\pi_t$ takes now the special form

$$\pi_t = \sum_{n=0}^{\mathcal{N}_t} \mathbb{1}_{\{t = \tau_n\}} \bar{\pi}_n = \sum_{n=0}^{\infty} \mathbb{1}_{\{t = \tau_n\}} \bar{\pi}_n,$$

where $\bar{\pi}_n$ is $\mathcal{F}_{\tau_n}$-measurable. As we can see, the resulting $\pi_t$ is nonzero only if we have an occurrence at $t$.

By decomposing $\bar{\pi}_n = \bar{\varpi}_n \bar{p}_n$ with $\sum_{n=0}^{\infty} \bar{\varpi}_n = 1$ and $\mathbb{E}_{\infty}[\bar{p}_n] = 1$, we can define the equivalent of all performance measures introduced in Section 3.1. We limit our presentation to Lorden's measure since this is the case we are going to treat in detail. If we use the last equation for $\pi_t$ in (3.1), we obtain the following form for our performance measure $\mathcal{J}(T)$:

$$(4.1) \qquad \mathcal{J}(T) = \frac{\sum_{n=0}^{\infty} \bar{\varpi}_n \mathbb{E}_\infty [\bar{p}_n \mathbb{E}_0 [(T - \tau_n)^+ | \mathcal{F}_{\tau_n}]]}{\sum_{n=0}^{\infty} \bar{\varpi}_n \mathbb{E}_\infty [\bar{p}_n \mathbb{1}_{\{T > \tau_n\}}]}.$$

Assuming no prior knowledge for $\{\bar{\varpi}_n\}$ and $\{\bar{p}_n\}$, we have to maximize $\mathcal{J}(T)$ in (4.1) with respect to the two processes. This leads to the following *extended* Lorden measure:

$$(4.2) \qquad \mathcal{J}_{\mathrm{EL}}(T) = \sup_{\{\bar{\varpi}_n\}, \{\bar{p}_n\}} \mathcal{J}(T) = \sup_{0 \le n} \mathrm{ess\,sup}\, \mathbb{E}_{\tau_n}[(T - \tau_n)^+ | \mathcal{F}_{\tau_n}].$$

The difference with the previous definition of $\mathcal{J}_{\mathrm{EL}}(T)$ in (3.10), is that the time instants $\tau_n$ are now s.t. instead of deterministic times.

### 4.1. *Detection of a change in the constant drift of a Brownian motion.*
Although it is possible to analyze the problem of detecting a change in the pdf of i.i.d. observations, we prefer to consider the continuous time alternative of detecting a change in the constant drift of a BM. This is because the corresponding solution is more elegant, offering formulas for the optimum performance and therefore allowing for direct comparison with the classical CUSUM test. Thus, let us assume that the observation process $\{\xi_t\}$ is a BM satisfying $\xi_t = \mu(t - \tau)^+ + w_t$, where $w_t$ a standard Wiener process and $\mu$ a known constant drift. For the change-time $\tau$ we assume that it can be equal to any $\tau_n$ from the observable sequence of s.t. $\{\tau_n\}$. Finally for the occurrence times $\{\tau_n\}$, we assume that they are Poisson distributed with a constant rate $\lambda$ and independent from the observation process $\{\xi_t\}$.

We recall that the problem of detecting a change in the drift of a BM has been considered with the classical Lorden measure (where occurrences are not taken into account, therefore the change is assumed to happen at any time instant) by Shiryaev (1996) and Beibel (1996) and under a more general framework by Moustakides (2004).

If we denote with $u_t$ the log-likelihood ratio between the two probability measures, then $u_t = -0.5\mu^2 t + \mu \xi_t$. Let us consider the following process $\{m_t\}$:

$$m_0 = 0; \qquad m_t = m_0 \wedge \left( \inf_{1 \le n \le \mathcal{N}_t} u_{\tau_n} \right) = \left( \inf_{1 \le n \le \mathcal{N}_t} u_{\tau_n} \right)^-,$$

where $x^- = \min\{x, 0\}$. Notice that $m_t$ starts from 0 and becomes the running minimum of the process $\{u_t\}$ but updated *only at the occurrence times*. We can now define the extended CUSUM (ECUSUM) process as follows:

$$y_t = u_t - m_t$$

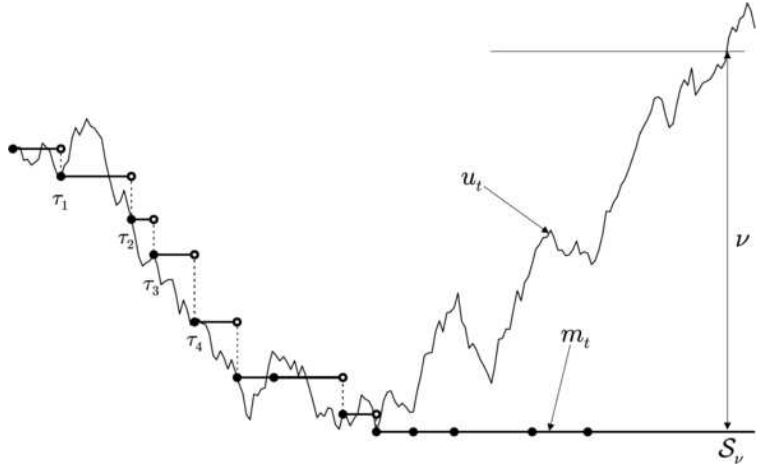

Fig. 1. *Sample paths of $u_t$, $m_t$ and stopping of ECUSUM.*

and the corresponding ECUSUM s.t. with threshold $\nu \geq 0$ as

$$\mathcal{S}_\nu = \inf_{0 < t}\{t : y_t \geq \nu\}.$$

As opposed to the CUSUM process which is always nonnegative, the ECUSUM process $y_t$ can take upon negative values as well.

Figure 1 depicts an example of the paths of $\{u_t\}$ and $\{m_t\}$. Process $\{m_t\}$ is piece-wise constant with right continuous paths and can exhibit jumps at the occurrence instances $\{\tau_n\}$. The ECUSUM process $\{y_t\}$, being the difference of $\{u_t\}$ (which is continuous) and $\{m_t\}$, is also right continuous with continuous paths between occurrences. From Figure 1, we can also deduce that $\{y_t\}$ exhibits a jump at $\tau_n$ only if $y_{\tau_{n-}} < 0$ in which case $y_{\tau_n}$ becomes 0. This can be written as

(4.3) $$y_{\tau_n} = (y_{\tau_{n-}})^+.$$

For technical reason, it is also necessary to introduce a version of ECUSUM which can start from any value $y_0 = y$ (as compared to the regular version which starts at $y_0 = 0$). For this we simply have to assume that $u_t = y - 0.5\mu^2 t + \mu\xi_t$, while $m_t, y_t$ and the s.t. are defined as before. To distinguish this new version of the s.t. from the regular one let us denote it as $\tilde{\mathcal{S}}_\nu$. It is then clear that $\mathcal{S}_\nu = \tilde{\mathcal{S}}_\nu$ for $y = 0$. Since inter-occurrence times are i.i.d. and independent from the past, the process $\{y_{\tau_n}\}$ is Markov and $\tilde{\mathcal{S}}_\nu$ given that $y_0 = y$ has the same statistics as $(\tilde{\mathcal{S}}_\nu - \tau_n)^+$ given that $y_{\tau_n} = y$ and $\tilde{\mathcal{S}}_\nu > \tau_n$.

4.2. *Performance evaluation of the ECUSUM test.* In this subsection, we are going to obtain a formula for the expectation of $\tilde{\mathcal{S}}_\nu$. We first present a lemma that states an important property for this quantity.

Lemma 3. *Let the occurrences be Poisson distributed with rate $\lambda$, then the average $\mathbb{E}[\tilde{\mathcal{S}}_\nu]$ is decreasing in $y_0 = y$ and for every $y$ we have*

$$(4.4) \qquad \mathbb{E}[\tilde{\mathcal{S}}_\nu | y_0 = y] \leq \frac{1}{\lambda} + \mathbb{E}[\mathcal{S}_\nu].$$

Proof. The paths of $\{y_t\}$ are increasing in $y_0 = y$, therefore $\tilde{\mathcal{S}}_\nu$ is decreasing in $y$ and so is $\mathbb{E}[\tilde{\mathcal{S}}_\nu]$; consequently for $y \geq 0$ we have $\mathbb{E}[\tilde{\mathcal{S}}_\nu | y_0 = y] \leq \mathbb{E}[\tilde{\mathcal{S}}_\nu | y_0 = 0] = \mathbb{E}[\mathcal{S}_\nu]$. Assume now that $y < 0$, then since $\tilde{\mathcal{S}}_\nu \leq \tau_1 + (\tilde{\mathcal{S}}_\nu - \tau_1)^+$, by taking expectation we can write

$$\begin{aligned}
\mathbb{E}[\tilde{\mathcal{S}}_\nu | y_0 = y] &\leq \mathbb{E}[\tau_1] + \mathbb{E}[(\tilde{\mathcal{S}}_\nu - \tau_1)^+ | y_0 = y] \\
&= \mathbb{E}[\tau_1] + \mathbb{E}[(\tilde{\mathcal{S}}_\nu - \tau_1)^+ | y_0 = y, \tilde{\mathcal{S}}_\nu > \tau_1] \mathbb{P}[\tilde{\mathcal{S}}_\nu > \tau_1 | y_0 = y] \\
&\leq \mathbb{E}[\tau_1] + \mathbb{E}[(\tilde{\mathcal{S}}_\nu - \tau_1)^+ | y_0 = y, \tilde{\mathcal{S}}_\nu > \tau_1] \\
&= \mathbb{E}[\tau_1] + \mathbb{E}[\mathbb{E}[(\tilde{\mathcal{S}}_\nu - \tau_1)^+ | y_{\tau_1}, \tilde{\mathcal{S}}_\nu > \tau_1] | y_0 = y] \\
&\leq \mathbb{E}[\tau_1] + \sup_{z \geq 0} \mathbb{E}[(\tilde{\mathcal{S}}_\nu - \tau_1)^+ | y_{\tau_1} = z, \tilde{\mathcal{S}}_\nu > \tau_1] \\
&= \frac{1}{\lambda} + \sup_{z \geq 0} \mathbb{E}[\tilde{\mathcal{S}}_\nu | y_0 = z] \\
&= \frac{1}{\lambda} + \mathbb{E}[\tilde{\mathcal{S}}_\nu | y_0 = 0] \\
&= \frac{1}{\lambda} + \mathbb{E}[\mathcal{S}_\nu].
\end{aligned}$$

Where we have used the property that $(\tilde{\mathcal{S}}_\nu - \tau_1)^+$ conditioned on the event that $\tilde{\mathcal{S}}_\nu > \tau_1$ and $y_{\tau_1} = y$, has the same statistics as $\tilde{\mathcal{S}}_\nu$ given $y_0 = y$ and, furthermore, that at an occurrence the ECUSUM statistics is nonnegative. This concludes the proof. □

From Lemma 3 we deduce that $\mathbb{E}[\tilde{\mathcal{S}}_\nu | y_0 = y]$ is decreasing and *uniformly bounded* in $y$. Let us now proceed with the computation of $\mathbb{E}[\tilde{\mathcal{S}}_\nu | y_0 = y]$. We have the following theorem that provides the desired formula.

Theorem 1. *Let $u_t = y + at + bw_t$ with $\{w_t\}$ a standard Wiener process with $w_0 = 0$ and $a, b \neq 0$. Define the ECUSUM s.t. $\tilde{\mathcal{S}}_\nu$ as above, with the occurrences being Poisson distributed with rate $\lambda$, then for $y \leq \nu$ the expectation of $\tilde{\mathcal{S}}_\nu$ is given by the following expression:*

$$(4.5) \qquad \mathbb{E}[\tilde{\mathcal{S}}_\nu | y_0 = y] = \begin{cases} \dfrac{1}{a}[-y + \nu + A(e^{-2ya/b^2} - e^{-2\nu a/b^2})], & y \geq 0, \\[2mm] \dfrac{1}{a}[\nu + A(1 - e^{-2\nu a/b^2})] + \dfrac{1}{\lambda}[1 - e^{ry}], & y < 0, \end{cases}$$

*where*

$$r = \frac{-a + \sqrt{a^2 + 2\lambda b^2}}{b^2}, \qquad A = \frac{b^2}{2a}\left(\frac{ar}{\lambda} - 1\right).$$

PROOF. Denote with $f(y)$ the function in the right-hand side of (4.5) which, as we can verify, is twice continuously differentiable, strictly decreasing in $y$ and uniformly bounded for $-\infty < y \leq \nu$. Consider now the difference $f(y_t) - f(y_0)$, we can then write

$$f(y_t) - f(y_0) = f(y_t) - f(y_{\tau_{\mathcal{N}_t}}) + \sum_{n=1}^{\mathcal{N}_t}[f(y_{\tau_n}) - f(y_{\tau_{n-1}})]$$

$$= f(y_t) - f(y_{\tau_{\mathcal{N}_t}}) + \sum_{n=1}^{\mathcal{N}_t}[f(y_{\tau_{n-}}) - f(y_{\tau_{n-1}})]$$

$$+ \sum_{n=1}^{\mathcal{N}_t}[f(y_{\tau_n}) - f(y_{\tau_{n-}})],$$

where we used the fact that $\{y_t\}$ is right continuous. In the time interval $[\tau_{n-1}, \tau_n)$, the process $y_t$ has continuous paths and $m_t$ is constant, therefore using Itô calculus we can write

$$f(y_{\tau_{n-}}) - f(y_{\tau_{n-1}}) = \int_{\tau_{n-1}}^{\tau_n} f'(y_{s-})(a\,ds + b\,dw_s) + 0.5b^2 f''(y_{s-})\,ds.$$

If $t$ is not an occurrence, a similar expression holds for the time interval $[\tau_{\mathcal{N}_t}, t]$. This suggests that

$$f(y_t) - f(y_{\tau_{\mathcal{N}_t}}) + \sum_{n=1}^{\mathcal{N}_t}[f(y_{\tau_{n-}}) - f(y_{\tau_{n-1}})]$$

$$= \int_0^t [af'(y_{s-}) + 0.5b^2 f''(y_{s-})]\,ds + \int_0^t bf'(y_{s-})\,dw_s.$$

The sum involving the jumps, using (4.3), can be written as

$$\sum_{n=1}^{\mathcal{N}_t}[f(y_{\tau_n}) - f(y_{\tau_{n-}})] = \sum_{n=1}^{\mathcal{N}_t}[f((y_{\tau_{n-}})^+) - f(y_{\tau_{n-}})]$$

$$= \int_0^t [f((y_{s-})^+) - f(y_{s-})]\,d\mathcal{N}_s.$$

Combining the two expressions leads to

$$f(y_t) - f(y_0) = \int_0^t [af'(y_{s-}) + 0.5b^2 f''(y_{s-})]\,ds + \int_0^t bf'(y_{s-})\,dw_s$$

(4.6)

$$+ \int_0^t [f((y_{s-})^+) - f(y_{s-})]\,d\mathcal{N}_s.$$

For any integer $n$ let $\tilde{\mathcal{S}}_\nu^n = \tilde{\mathcal{S}}_\nu \wedge n$. Then we know that for a process $\{\omega_t\}$ which is a $\{\mathcal{F}_t\}$-adapted and uniformly bounded in the sense that $|\omega_t| \le c < \infty$, we have from Protter (2004) that $\mathbb{E}[\int_0^{\tilde{\mathcal{S}}_\nu^n} \omega_{s-}\, d\mathcal{N}_s] = \mathbb{E}[\int_0^{\tilde{\mathcal{S}}_\nu^n} \omega_{s-}\lambda\, ds]$ and from Karatzas and Shreve (1988) that $\mathbb{E}[\int_0^{\tilde{\mathcal{S}}_\nu^n} \omega_{s-}\, dw_s] = 0$. Replacing $t$ with the s.t. $\tilde{\mathcal{S}}_\nu^n$ in (4.6), taking expectation and using the fact that $f(y^+) - f(y)$ and $f'(y)$ are uniformly bounded for $y \in (-\infty, \nu]$, allows us to write

$$\mathbb{E}[f(y_{\tilde{\mathcal{S}}_\nu^n})] - f(y_0)$$
$$= \mathbb{E}\left[\int_0^{\tilde{\mathcal{S}}_\nu^n} \{af'(y_{t-}) + 0.5b^2 f''(y_{t-}) + \lambda[f((y_{t-})^+) - f(y_{t-})]\}\, dt\right].$$

It is straightforward to verify that the function $f(y)$ is a solution to the differential equation

$$(4.7) \quad af'(y) + 0.5b^2 f''(y) + \lambda[f(y^+) - f(y)] = -1, \qquad -\infty < y \le \nu.$$

This, if substituted in the previous expression, yields

$$f(y_0) - \mathbb{E}[f(y_{\tilde{\mathcal{S}}_\nu^n})] = \mathbb{E}[\tilde{\mathcal{S}}_\nu^n].$$

Letting now $n \to \infty$, we have $\tilde{\mathcal{S}}_\nu^n \to \tilde{\mathcal{S}}_\nu$ monotonically. In the previous equality, using monotone convergence on the right-hand side and bounded convergence [since $f(y)$ is uniformly bounded] on the left, we obtain

$$f(y_0) - \mathbb{E}[f(y_{\tilde{\mathcal{S}}_\nu})] = \mathbb{E}[\tilde{\mathcal{S}}_\nu].$$

At the time of stopping the process $\{y_t\}$ hits the threshold $\nu$ (see Figure 1), therefore, we have $y_{\tilde{\mathcal{S}}_\nu} = \nu$, suggesting that $\mathbb{E}[f(y_{\tilde{\mathcal{S}}_\nu})] = f(\nu)$. We can now verify that $f(\nu) = 0$, which yields $f(y_0) = \mathbb{E}[\tilde{\mathcal{S}}_\nu]$ and completes the proof. $\square$

REMARK 1. One might wonder, why is (4.5) the desired formula and not any other solution of the differential equation in (4.7) that satisfies the boundary condition $f(\nu) = 0$? It turns out that among the solutions of (4.7) that are twice continuously differentiable (property needed to apply Itô calculus) and satisfy the boundary condition $f(\nu) = 0$, the formula in (4.5) is the *unique* solution which is uniformly bounded in $(-\infty, \nu]$ (property imposed by Lemma 3).

By letting $\lambda \to \infty$ and setting $y = 0$ in (4.5), we recover the average run length of the classical CUSUM test as obtained in Taylor (1975). If we denote by $g_\nu(y), h_\nu(y)$ the average of $\tilde{\mathcal{S}}_\nu$ under $\mathbb{P}_0$ and $\mathbb{P}_\infty$ respectively, then under $\mathbb{P}_0$ we have $u_t = y - 0.5\mu^2 t + \mu \xi_t = y + 0.5\mu^2 t + \mu w_t$, therefore by substituting $a = 0.5\mu^2, b = \mu$ in (4.5), we can write

$$g_\nu(y) = \mathbb{E}_0[\tilde{\mathcal{S}}_\nu | y_0 = y] = \begin{cases} \dfrac{2}{\mu^2}[-y + \nu + A_0(e^{-y} - e^{-\nu})], & y \geq 0, \\[2mm] \dfrac{2}{\mu^2}[\nu + A_0(1 - e^{-\nu})] + \dfrac{1}{\lambda}[1 - e^{r_0 y}], & y < 0, \end{cases}$$

where

$$r_0 = -\frac{1}{2} + \sqrt{\frac{1}{4} + \frac{2\lambda}{\mu^2}}, \qquad A_0 = \frac{\mu^2}{2\lambda} r_0 - 1.$$

Similarly substituting $a = -0.5\mu^2, b = \mu$ in (4.5), we obtain

$$h_\nu(y) = \mathbb{E}_\infty[\tilde{\mathcal{S}}_\nu | y_0 = y] = \begin{cases} \dfrac{2}{\mu^2}[y - \nu + A_\infty(e^\nu - e^y)], & y \geq 0, \\[2mm] \dfrac{2}{\mu^2}[-\nu + A_\infty(e^\nu - 1)] + \dfrac{1}{\lambda}[1 - e^{r_\infty y}], & y < 0, \end{cases}$$

where

$$r_\infty = \frac{1}{2} + \sqrt{\frac{1}{4} + \frac{2\lambda}{\mu^2}}, \qquad A_\infty = \frac{\mu^2}{2\lambda} r_\infty + 1.$$

To compute the performance of the regular ECUSUM s.t. $\mathcal{S}_\nu$ we must set $y = 0$ in the previous formulas. It is then clear that $g_\nu(0)$ expresses the (worst) average detection delay and $h_\nu(0)$ the average period between false alarms for $\mathcal{S}_\nu$. Specifically, after noticing that $r_0 r_\infty = 2\lambda/\mu^2$, we have

$$(4.8) \qquad g_\nu(0) = \mathbb{E}_0[\mathcal{S}_\nu] = \frac{2}{\mu^2}\left\{[\nu - 1 + e^{-\nu}] + \frac{1}{r_\infty}(1 - e^{-\nu})\right\},$$

$$(4.9) \qquad h_\nu(0) = \mathbb{E}_\infty[\mathcal{S}_\nu] = \frac{2}{\mu^2}\left\{[e^\nu - \nu - 1] + \frac{1}{r_0}(e^\nu - 1)\right\},$$

where the first term in both right hand side expressions corresponds to the performance of the classical CUSUM test (obtained by letting $\lambda \to \infty$).

Figure 2 depicts the normalized average detection delay $\mu^2 g_\nu(0)/2$ as a function of the normalized average false alarm period $\mu^2 h_\nu(0)/2$, for different values of the ratio $\mu^2/\lambda$. We observe that, in the average, ECUSUM detects the change faster than CUSUM. Of course this is not surprising since ECUSUM has available more information than CUSUM (CUSUM does not observe the occurrences). We can also see that the performance difference between the two schemes, for given value of $\mu^2/\lambda$, is uniformly bounded by a constant. Finally, we conclude that the gain obtained by using ECUSUM instead of CUSUM becomes significant only for large values of the parameter $\mu^2/\lambda$ or, equivalently, when the occurrences that can trigger the change are very *infrequent*.

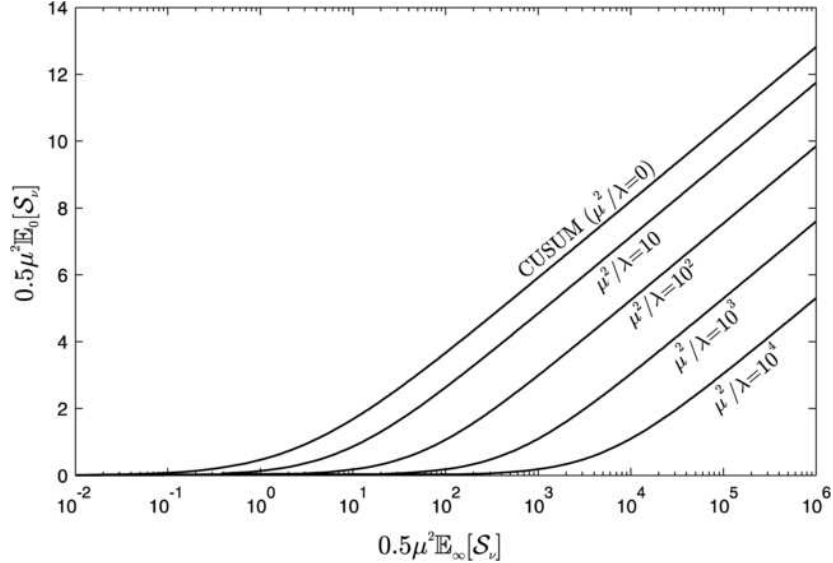

FIG. 2.  *Normalized average detection delay of ECUSUM as a function of normalized average false alarm period, for different values of the parameter $\mu^2/\lambda$.*

**5. Optimality of ECUSUM.** Using the formula in (4.9) for the average period between false alarms, we can relate the threshold $\nu$ to the false alarm constraint parameter $\gamma$ through the equation

$$h_\nu(0) = \frac{2}{\mu^2}\left\{[e^\nu - \nu - 1] + \frac{1}{r_0}(e^\nu - 1)\right\} = \gamma.$$

The left-hand side of the last equality is increasing in $\nu$ and for $\nu = 0$ it is equal to 0, also for $\nu \to \infty$ it tends to $\infty$, we can, therefore, conclude that for given $\gamma \geq 0$, the last equation has a unique solution which we call $\nu_\star$. Since $\nu_\star$ is the solution to the previous equation it is clear that

$$(5.1) \qquad\qquad\qquad h_{\nu_\star}(0) = \gamma.$$

Using $\nu_\star$ as threshold we can now define the corresponding ECUSUM s.t. $\mathcal{S}_{\nu_\star}$. Our goal in the sequel is to demonstrate that this test is optimum in the extended Lorden sense. We recall that the occurrence times are Poisson distributed with a *known* constant rate $\lambda$. Observe however that $\lambda$ enters only in the correspondence between the threshold $\nu_\star$ and the constraint $\gamma$ *without affecting* the ECUSUM test otherwise.

Consider the functions $g_{\nu_\star}(y), h_{\nu_\star}(y)$ associated with $\mathcal{S}_{\nu_\star}$. Both functions will play a key role in our proof of optimality. The next lemma presents an important property for each function which is an immediate consequence of Theorem 1.

LEMMA 4. *If $T$ is a s.t. and $\mathcal{S}_\nu$ the regular ECUSUM s.t. with threshold $\nu$, define $T_\nu = T \wedge \mathcal{S}_\nu$, then*

$$(5.2) \qquad \mathbb{E}_\infty[h_{\nu_\star}(0) - h_{\nu_\star}(y_{T_\nu})] = \mathbb{E}_\infty[T_\nu],$$

$$(5.3) \quad \mathbb{E}_{\tau_n}[\{g_{\nu_\star}(y_{\tau_n}) - g_{\nu_\star}(y_{T_\nu})\}\mathbb{1}_{\{T_\nu > \tau_n\}}|\mathcal{F}_{\tau_n}] = \mathbb{E}_{\tau_n}[(T_\nu - \tau_n)^+|\mathcal{F}_{\tau_n}].$$

With the next theorem we provide a suitable lower bound for the extended Lorden measure. First, we introduce a technical lemma.

LEMMA 5. *Let $T$ be a s.t. and define $T_\nu = T \wedge \mathcal{S}_\nu$, then*

$$\mathbb{E}_\infty[e^{y_{T_\nu}}] \geq 1.$$

PROOF. Following similar steps as in the proof for Theorem 1, if $f(y)$ is a twice continuously differentiable function with $f(y^+) - f(y)$ and $f'(y)$ uniformly bounded for $y \leq \nu$, we have

$$\mathbb{E}_\infty[f(y_{T_\nu})] - f(y_0)$$
$$= \mathbb{E}_\infty\left[\int_0^{T_\nu}\{-0.5\mu^2 f'(y_{t-}) + 0.5\mu^2 f''(y_{t-}) + \lambda[f((y_{t-})^+) - f(y_{t-})]\}\,dt\right].$$

For $f(y) = e^y$ and recalling that we treat the regular ECUSUM case with $y_0 = 0$, we immediately obtain that

$$\mathbb{E}_\infty[e^{y_{T_\nu}}] - 1 = \mathbb{E}_\infty\left[\int_0^{T_\nu}\lambda(e^{y_{s-}^+} - e^{y_{s-}})\,ds\right]$$
$$= \lambda\mathbb{E}_\infty\left[\int_0^{T_\nu}(1 - e^{y_{s-}})^+\,ds\right] \geq 0.$$

This concludes the proof. □

THEOREM 2. *For any s.t. $T$ let $T_\nu = T \wedge \mathcal{S}_\nu$, then*

$$(5.4) \qquad \mathcal{J}_{\mathrm{EL}}(T) \geq g_{\nu_\star}(0) - \frac{\mathbb{E}_\infty[g_{\nu_\star}(y_{T_\nu})e^{y_{T_\nu}}]}{\mathbb{E}_\infty[e^{y_{T_\nu}}]}.$$

PROOF. The proof follows similar steps as in Theorem 2, Moustakides (2004). Since $T \geq T_\nu$, it is clear that $\mathcal{J}_{\mathrm{EL}}(T) \geq \mathcal{J}_{\mathrm{EL}}(T_\nu)$. Also from the definition of $\mathcal{J}_{\mathrm{EL}}(\cdot)$ in (4.2) we conclude

$$(5.5) \quad \mathcal{J}_{\mathrm{EL}}(T) \geq \mathcal{J}_{\mathrm{EL}}(T_\nu) \geq \mathbb{E}_{\tau_n}[(T_\nu - \tau_n)^+|\mathcal{F}_{\tau_n}], \qquad n = 0, 1, 2, \ldots.$$

Using (5.3) from Lemma 4, the previous inequality for $n \geq 1$ can be written as

$$\mathcal{J}_{\mathrm{EL}}(T) \geq \mathbb{E}_{\tau_n}[(T_\nu - \tau_n)^+|\mathcal{F}_{\tau_n}]$$
$$= \mathbb{E}_{\tau_n}[\{g_{\nu_\star}(y_{\tau_n}) - g_{\nu_\star}(y_{T_\nu})\}\mathbb{1}_{\{T_\nu > \tau_n\}}|\mathcal{F}_{\tau_n}]$$
$$= \mathbb{E}_\infty[e^{u_{T_\nu} - u_{\tau_n}}\{g_{\nu_\star}(y_{\tau_n}) - g_{\nu_\star}(y_{T_\nu})\}\mathbb{1}_{\{T_\nu > \tau_n\}}|\mathcal{F}_{\tau_n}].$$

Multiplying both sides with the nonnegative quantity $(1 - e^{m_{\tau_n} - m_{\tau_{n-1}}})\mathbb{1}_{\{T_\nu > \tau_n\}}$ and taking expectation with respect to $\mathbb{P}_\infty$ yields

$$
\begin{aligned}
(5.6) \quad & \mathcal{J}_{\mathrm{EL}}(T)\mathbb{E}_\infty[(1 - e^{m_{\tau_n} - m_{\tau_{n-1}}})\mathbb{1}_{\{T_\nu > \tau_n\}}] \\
& \geq \mathbb{E}_\infty[e^{u_{T_\nu} - u_{\tau_n}}(1 - e^{m_{\tau_n} - m_{\tau_{n-1}}})\{g_{\nu_\star}(y_{\tau_n}) - g_{\nu_\star}(y_{T_\nu})\}\mathbb{1}_{\{T_\nu > \tau_n\}}].
\end{aligned}
$$

Now notice that $1 - e^{m_{\tau_n} - m_{\tau_{n-1}}}$ is different from 0, only when there is a jump in $m_t$ at $\tau_n$, in which case $y_{\tau_n} = 0$ and $m_{\tau_n} = u_{\tau_n}$. This means that

$$
\begin{aligned}
& \mathbb{E}_\infty[e^{u_{T_\nu} - u_{\tau_n}}(1 - e^{m_{\tau_n} - m_{\tau_{n-1}}})\{g_{\nu_\star}(y_{\tau_n}) - g_{\nu_\star}(y_{T_\nu})\}\mathbb{1}_{\{T_\nu > \tau_n\}}] \\
& = \mathbb{E}_\infty[e^{u_{T_\nu} - m_{\tau_n}}(1 - e^{m_{\tau_n} - m_{\tau_{n-1}}})\{g_{\nu_\star}(0) - g_{\nu_\star}(y_{T_\nu})\}\mathbb{1}_{\{T_\nu > \tau_n\}}] \\
& = \mathbb{E}_\infty[e^{u_{T_\nu}}(e^{-m_{\tau_n}} - e^{-m_{\tau_{n-1}}})\{g_{\nu_\star}(0) - g_{\nu_\star}(y_{T_\nu})\}\mathbb{1}_{\{T_\nu > \tau_n\}}].
\end{aligned}
$$

Furthermore, since $\mathbb{E}_\infty[e^{u_{T_\nu} - u_{\tau_n}}|\mathcal{F}_{\tau_n}] = 1$ we can write

$$
\begin{aligned}
\mathbb{E}_\infty[(1 - e^{m_{\tau_n} - m_{\tau_{n-1}}})\mathbb{1}_{\{T_\nu > \tau_n\}}] & = \mathbb{E}_\infty[e^{u_{T_\nu} - u_{\tau_n}}(1 - e^{m_{\tau_n} - m_{\tau_{n-1}}})\mathbb{1}_{\{T_\nu > \tau_n\}}] \\
& = \mathbb{E}_\infty[e^{u_{T_\nu} - m_{\tau_n}}(1 - e^{m_{\tau_n} - m_{\tau_{n-1}}})\mathbb{1}_{\{T_\nu > \tau_n\}}] \\
& = \mathbb{E}_\infty[e^{u_{T_\nu}}(e^{-m_{\tau_n}} - e^{-m_{\tau_{n-1}}})\mathbb{1}_{\{T_\nu > \tau_n\}}].
\end{aligned}
$$

Substituting the two equalities in (5.6) and summing over all $n \geq 1$ we have

$$
\begin{aligned}
(5.7) \quad & \mathcal{J}_{\mathrm{EL}}(T)\sum_{n=1}^\infty \mathbb{E}_\infty[e^{u_{T_\nu}}(e^{-m_{\tau_n}} - e^{-m_{\tau_{n-1}}})\mathbb{1}_{\{T_\nu > \tau_n\}}] \\
& \geq \sum_{n=1}^\infty \mathbb{E}_\infty[e^{u_{T_\nu}}(e^{-m_{\tau_n}} - e^{-m_{\tau_{n-1}}})\{g_{\nu_\star}(0) - g_{\nu_\star}(y_{T_\nu})\}\mathbb{1}_{\{T_\nu > \tau_n\}}].
\end{aligned}
$$

In the second sum, interchanging summation and expectation, yields

$$
\begin{aligned}
& \sum_{n=1}^\infty \mathbb{E}_\infty[e^{u_{T_\nu}}(e^{-m_{\tau_n}} - e^{-m_{\tau_{n-1}}})\{g_{\nu_\star}(0) - g_{\nu_\star}(y_{T_\nu})\}\mathbb{1}_{\{T_\nu > \tau_n\}}] \\
& = \mathbb{E}_\infty\left[\{g_{\nu_\star}(0) - g_{\nu_\star}(y_{T_\nu})\}e^{u_{T_\nu}}\sum_{n=1}^{\mathcal{N}_{T_\nu}}(e^{-m_{\tau_n}} - e^{-m_{\tau_{n-1}}})\right] \\
& = \mathbb{E}_\infty[\{g_{\nu_\star}(0) - g_{\nu_\star}(y_{T_\nu})\}e^{u_{T_\nu}}(e^{-m_{T_\nu}} - 1)].
\end{aligned}
$$

Similarly for the first sum we have

$$
\begin{aligned}
& \sum_{n=1}^\infty \mathbb{E}_\infty[e^{u_{T_\nu}}(e^{-m_{\tau_n}} - e^{-m_{\tau_{n-1}}})\mathbb{1}_{\{T_\nu > \tau_n\}}] \\
& = \mathbb{E}_\infty[e^{u_{T_\nu}}(e^{-m_{T_\nu}} - 1)].
\end{aligned}
$$

Substituting the two expressions in (5.7) we obtain

$$
\begin{aligned}
(5.8) \quad & \mathcal{J}_{\mathrm{EL}}(T)\mathbb{E}_{\infty}[e^{u_{T_{\nu}}}(e^{-m_{T_{\nu}}}-1)] \\
& \geq \mathbb{E}_{\infty}[\{g_{\nu_{\star}}(0)-g_{\nu_{\star}}(y_{T_{\nu}})\}e^{u_{T_{\nu}}}(e^{-m_{T_{\nu}}}-1)].
\end{aligned}
$$

There is one last inequality we have not used so far from (5.5), namely for $n=0$. Recalling that $\tau_0=0$, this inequality takes the form

$$
(5.9) \qquad\qquad \mathcal{J}_{\mathrm{EL}}(T) \geq \mathbb{E}_0[T_{\nu}].
$$

Using (5.3) from Lemma 4, we get

$$
\begin{aligned}
\mathbb{E}_0[T_{\nu}] &= \mathbb{E}_0[g_{\nu_{\star}}(0)-g_{\nu_{\star}}(y_{T_{\nu}})] \\
&= \mathbb{E}_{\infty}[\{g_{\nu_{\star}}(0)-g_{\nu_{\star}}(y_{T_{\nu}})\}e^{u_{T_{\nu}}}].
\end{aligned}
$$

Also since $\mathbb{E}_{\infty}[e^{u_{T_{\nu}}}]=1$, (5.9) is equivalent to

$$
\mathcal{J}_{\mathrm{EL}}(T)\mathbb{E}_{\infty}[e^{u_{T_{\nu}}}] \geq \mathbb{E}_{\infty}[\{g_{\nu_{\star}}(0)-g_{\nu_{\star}}(y_{T_{\nu}})\}e^{u_{T_{\nu}}}].
$$

If this is added to (5.8), we obtain

$$
\mathcal{J}_{\mathrm{EL}}(T)\mathbb{E}_{\infty}[e^{u_{T_{\nu}}-m_{T_{\nu}}}] \geq \mathbb{E}_{\infty}[\{g_{\nu_{\star}}(0)-g_{\nu_{\star}}(y_{T_{\nu}})\}e^{u_{T_{\nu}}-m_{T_{\nu}}}],
$$

or

$$
(5.10) \quad \mathcal{J}_{\mathrm{EL}}(T)\mathbb{E}_{\infty}[e^{y_{T_{\nu}}}] \geq g_{\nu_{\star}}(0)\mathbb{E}_{\infty}[e^{y_{T_{\nu}}}]-\mathbb{E}_{\infty}[g_{\nu_{\star}}(y_{T_{\nu}})e^{y_{T_{\nu}}}].
$$

Since $y_{T_{\nu}} \leq \nu$, thanks also to Lemma 5, we have $e^{\nu} \geq \mathbb{E}_{\infty}[e^{y_{T_{\nu}}}] \geq 1$. We can thus divide both sides of (5.10) with $\mathbb{E}_{\infty}[e^{y_{T_{\nu}}}]$ and obtain the desired expression. $\square$

We will base our proof of optimality of $\mathcal{S}_{\nu_{\star}}$ on Theorem 2. Let us first introduce an additional technical lemma.

LEMMA 6.    *If $T$ is a s.t. and $T_{\nu}=T\wedge\mathcal{S}_{\nu}$ then the function $\psi(\nu)=\mathbb{E}_{\infty}[T_{\nu}]$ is continuous and increasing in $\nu$ with $\psi(0)=0$ and $\psi(\infty)=\mathbb{E}_{\infty}[T]$.*

PROOF.    The proof is exactly similar as in Lemma 3, Moustakides (2004). Consider $\kappa > \nu \geq 0$, then

$$
0 \leq T_{\kappa}-T_{\nu} \leq \mathcal{S}_{\kappa}-\mathcal{S}_{\nu},
$$

from which we obtain

$$
0 \leq \psi(\kappa)-\psi(\nu) \leq h_{\kappa}(0)-h_{\nu}(0).
$$

From (4.9) we have that the function $h_{\nu}(0)$ is continuous in $\nu$. If we use this property in the previous relation, we deduce that $\psi(\nu)$ is also continuous.

Finally, we can directly verify that the two limiting values $\psi(0), \psi(\infty)$ are correct.   $\square$

An immediate consequence of the previous lemma is the fact that if $\mathbb{E}_\infty[T] > \gamma$ then we can find a threshold $\nu$ such that $\mathbb{E}_\infty[T_\nu] = \gamma$. Since $J_{\mathrm{EL}}(T) \geq J_{\mathrm{EL}}(T_\nu)$ this suggests that for the proof of optimality of $\mathcal{S}_{\nu_\star}$ we can limit ourselves to s.t. that satisfy the false alarm constraint with equality.

THEOREM 3. *If a s.t. $T$ satisfies $\mathbb{E}_\infty[T] = \gamma$ then it possesses an extended Lorden measure $\mathcal{J}_{\mathrm{EL}}(T)$ that is no less than $g_{\nu_\star}(0) = \mathcal{J}_{\mathrm{EL}}(\mathcal{S}_{\nu_\star})$.*

PROOF. From $\mathbb{E}_\infty[T] = \gamma$, thanks to Lemma 6, for every $\epsilon > 0$ we can find a threshold $\nu_\epsilon$ such that for $T_{\nu_\epsilon} = T \wedge \mathcal{S}_{\nu_\epsilon}$ we have

$$(5.11) \qquad \gamma \geq \mathbb{E}_\infty[T_{\nu_\epsilon}] \geq \gamma - \epsilon.$$

Using (5.2) from Lemma 4, we can write

$$\mathbb{E}_\infty[h_{\nu_\star}(0) - h_{\nu_\star}(y_{T_{\nu_\epsilon}})] = \mathbb{E}_\infty[T_{\nu_\epsilon}].$$

Recalling from (5.1) that $h_{\nu_\star}(0) = \gamma$ and using (5.11) in the previous equality, we obtain

$$(5.12) \qquad \epsilon \geq \mathbb{E}_\infty[h_{\nu_\star}(y_{T_{\nu_\epsilon}})] \geq 0.$$

Define the function $p(y) = e^y g_{\nu_\star}(y) - (r_0/r_\infty)h_{\nu_\star}(y)$ and consider the derivative $p'(y)$. By direct substitution we can verify that $e^y g'_{\nu_\star}(y) = (r_0/r_\infty)h'_{\nu_\star}(y)$, from which we deduce that

$$p'(y) = e^y g_{\nu_\star}(y).$$

Since $g_{\nu_\star}(y)$ is strictly decreasing in $y$ and also satisfies $g_{\nu_\star}(\nu_\star) = 0$, this suggests that $p'(y)$ has the same sign as $\nu_\star - y$, or that $p(y)$ has a global maximum at $y = \nu_\star$. Because $p(\nu_\star) = 0$ this means that $p(y) \leq 0$ which yields $\mathbb{E}_\infty[p(y_{T_{\nu_\epsilon}})] \leq 0$. Using this inequality and replacing $p(y)$ by its definition, we obtain

$$0 \geq \mathbb{E}_\infty[p(y_{T_{\nu_\epsilon}})] = \mathbb{E}_\infty\left[e^{y_{T_{\nu_\epsilon}}} g_{\nu_\star}(y_{T_{\nu_\epsilon}}) - \frac{r_0}{r_\infty}h_{\nu_\star}(y_{T_{\nu_\epsilon}})\right]$$

$$\geq \mathbb{E}_\infty[e^{y_{T_{\nu_\epsilon}}} g_{\nu_\star}(y_{T_{\nu_\epsilon}})] - \frac{r_0}{r_\infty}\epsilon,$$

where for the last inequality we used (5.12). This yields

$$(5.13) \qquad \frac{r_0}{r_\infty}\epsilon \geq \mathbb{E}_\infty[e^{y_{T_{\nu_\epsilon}}} g_{\nu_\star}(y_{T_{\nu_\epsilon}})].$$

From Theorem 2, using (5.13) and Lemma 5, we can now write

$$\mathcal{J}_{\mathrm{EL}}(T) \geq g_{\nu_\star}(0) - \frac{\mathbb{E}_\infty[e^{y T_{\nu_\epsilon}} g_{\nu_\star}(y_{T_{\nu_\epsilon}})]}{\mathbb{E}_\infty[e^{y T_{\nu_\epsilon}}]} \geq g_{\nu_\star}(0) - \frac{r_0}{r_\infty}\epsilon.$$

Since $\epsilon$ is arbitrary this means that $\mathcal{J}_{\mathrm{EL}}(T) \geq g_{\nu_\star}(0)$, thus establishing optimality of ECUSUM. $\quad\square$

Department of Electrical
and Computer Engineering
University of Patras
26500 Rio
Greece
E-mail: moustaki@ece.upatras.gr